\documentclass[10pt,a4paper,fleqn]{article}
\usepackage{amsmath}
\usepackage{amssymb}
\usepackage{amscd}
\usepackage{plain}
\usepackage{graphicx}
\usepackage{epsfig}
\usepackage{fancyhdr}
\usepackage{t1enc,epsfig}
\usepackage[mathscr]{eucal}
\usepackage{epstopdf}

\newtheorem{theorem}{Theorem}

\newtheorem{corollary}[theorem]{Corollary}
\newtheorem{definition}[theorem]{Definition}
\newtheorem{example}[theorem]{Example}
\newtheorem{lemma}[theorem]{Lemma}
\newtheorem{proposition}[theorem]{Proposition}
\newtheorem{remark}[theorem]{Remark}
\newtheorem{conjecture}[theorem]{Conjecture}

\newenvironment{proof}{{\bf Proof. }}{\hfill$\rule{1ex}{1ex}$\par\medskip}
\begin{document}

\newcommand{\bt}{\begin{theorem}}
\newcommand{\et}{\end{theorem}}
\newcommand{\bd}{\begin{definition}}
\newcommand{\ed}{\end{definition}}
\newcommand{\bs}{\begin{proposition}}
\newcommand{\es}{\end{proposition}}
\newcommand{\bp}{\begin{proof}}
\newcommand{\ep}{\end{proof}}
\newcommand{\be}{\begin{equation}}
\newcommand{\ee}{\end{equation}}
\newcommand{\ul}{\underline}
\newcommand{\br}{\begin{remark}}
\newcommand{\er}{\end{remark}}
\newcommand{\bex}{\begin{example}}
\newcommand{\eex}{\end{example}}
\newcommand{\bc}{\begin{corollary}}
\newcommand{\ec}{\end{corollary}}
\newcommand{\bl}{\begin{lemma}}
\newcommand{\el}{\end{lemma}}
\newcommand{\bj}{\begin{conjecture}}
\newcommand{\ej}{\end{conjecture}}

\title{The Moufang's theorem for non-Moufang loops}

\author{Izabella Stuhl}

\date{}
\footnotetext{\par}

\maketitle

\begin{abstract}
We introduce a class of non-Moufang loops satisfying the Moufang's theorem.
\end{abstract}

\section{Introduction}
This note has two sources of motivation. First, it is a problem proposed by Rajah at the Loops '11 Conference. We say that a variety of loops $\mathcal{V}$ satisfies the Moufang's theorem if for every loop $L$ in $\mathcal{V}$ the following implication holds: for every elements $x, y, z \in L$ satisfying $x(yz) = (xy)z$ the subloop generated by $x, y, z$ is a group \cite{Mo}. Is it true that every variety satisfying the Moufang's theorem is contained in the variety of Moufang loops? The other source of motivation is an example (a Steiner loop of order $10$) presented in the talk of Giuliani at the 3rd Mile High Conference on Nonassociative Mathematics, details of which can be found in \cite{L}.

In \cite{CLRS} a large class of loops has been examined, formed by Steiner loops, and it was
shown that a Steiner loop satisfies the Moufang's theorem if and only if its corresponding Steiner triple system has
the property that every Pasch configuration generates a subsystem of order $7$.

In this paper we introduce a class of non-Moufang loops which is a class of non-commutative flexible alternative inverse property loops and show that satisfies the Moufang's theorem.

\section{Preliminaries and Result}

A {\it quasigroup} $(Q, \cdot)$ is a set $Q$ with a binary operation $(x,y)\mapsto x\cdot y$ such that
for given $a, b\in Q$ the equations $x\cdot a=b$ and $a\cdot y=b$ are uniquely solvable in $Q$. A quasigroup is called a {\it loop} if there is an identity element in $Q$. An idempotent totally symmetric quasigroup is called a {\it Steiner quasigroup}. A totally symmetric loop of exponent $2$ is called a {\it Steiner loop}. A {\it Moufang loop} is a loop satisfying any (and hence all) of the identities
$$
xy\cdot zx = x(yz\cdot x)\qquad (xy\cdot x)z = x(y\cdot xz)\qquad (zx\cdot y)x = z(x\cdot yx)\,.
$$
Moufang loops satisfy Moufang's theorem \cite{Mo}.

A {\it Steiner triple system} (STS) $\mathfrak{S}$ is a $2-(n,3,1)$ design, i.e., an incidence
structure consisting of points and blocks such that any two
distinct points are contained in precisely one block and any block
has precisely three points. A finite STS with
$n$ points exists if and only if $n \equiv 1$ or $3$ $(mod\hskip 4pt 6)$.

An STS $\mathfrak{S}$
generates a multiplication on pairs of different points $x, y$ taking as product the third point of the block joining $x$ and $y$.
Defining $x\cdot x=x$ we get a Steiner quasigroup associated with $\mathfrak{S}$. Further, adjoining an element $e$ with $ex=xe=x$, $xx=e$, yields a Steiner loop $S$.
Conversely, a Steiner loop $S$ determines an STS whose points are the elements of $S\setminus \{e\}$ and
the blocks are the triples $\{x, y, xy\}$ for all $x\neq y\in S\setminus \{e\}$.

A possible way to enrich STSs is to introduce in each block a cyclic order.
These objects occur in the literature as {\it oriented} STSs
$(\mathfrak{S}, T)$ (see \cite{W}). A method to relate quasigroups and loops to oriented STSs by applying Schreier-type extensions of quasigroups and loops (see \cite{NS2}, \cite{NSt}, \cite{NSt2}, \cite{S}) was proposed in \cite{S2}, \cite{SS2}. An {\it oriented Steiner loop} $L$ is a loop, for which there is an oriented STS $(\mathfrak{S},T)$
with the following properties (1)--(3).  (1) $L$ is a loop extension of the group of order $2$ by a Steiner loop $S=\mathfrak{S}\cup e$; (2) the restriction of the factor system of the extension to
$(\mathfrak{S}\times \mathfrak{S})\setminus \{(x,x), x\in
\mathfrak{S}\}$ coincides with the orientation function of
$(\mathfrak{S},T)$; and (3) $f(x,x)=-1$, respectively $f(x,x)=1$ $\forall$
$x\in S\setminus \{e\}$ (Cf. Definition 4 in \cite{SS}, p. 136.). With any oriented STS there
are associated an oriented Steiner loop of exponent $2$ and an oriented Steiner loop of exponent $4$.

We focus upon Steiner loops arising from special STSs,
namely, from Hall triple systems. A {\it Hall triple system} (HTS) is an STS in which every three non-collinear
points generate an affine plane over the field $GF(3)$; cf.  \cite{H}. These geometries are of interest
because they are the only known examples of perfect matroid designs other than classical projective and affine
geometries over finite fields and $t-(n,k,1)$ designs. The class of HTSs is an exceptional case to the Buekenhout's
characterization theorem (\cite{Bu} p. 368); they are the only non-degenerate non-affine "locally affine" geometries.
Furthermore, the "local affinity" of a non-affine geometry cannot be extended only as far as to dimension $3$. In
addition, HTSs are the only geometries of this type; see Theorem 1.3 in \cite{RR} p. 130.

In the present work, Steiner loops associated to HTSs are called {\it Hall loops} and oriented Steiner loops associated
to oriented HTSs are called {\it oriented Hall loops} (OHLs). Hall loops are not groups (cf. Theorem in \cite{DP} p. 250)
and not Moufang loops (cf \cite{Br}, Lemma 3.2, or Figure 1 below). Multiplication groups of such loops have been
determined in \cite{SS}.

\bt
Oriented Hall loops of exponent $4$ satisfy the Moufang's theorem.
\et

\bp
Let $L$ be an OHL, i.e., a loop extension of the group $Z_2$ of order $2$ by a Hall loop $S$, and let $\mathrm{x}=(x,\xi), \mathrm{y}=(y,\eta), \mathrm{z}=(z,\zeta)$ be elements of $L$. Consider
the associative law
$$[(x,\xi)\cdot (y,\eta)]\cdot (z,\zeta) = (x,\xi)\cdot [(y,\eta)\cdot (z,\zeta)].$$

First, we assume that $x\neq y\neq z\neq x$ and none among $x, y, z\in S$ is the product of the remaining
two elements, i.e., $x, y, z$ are non-collinear points of the corresponding HTS. Thus, $x, y, z$ generate an affine
plane over the field $GF(3)$. Hence, $x(yz)$ and $(xy)z$ are different points of the system (see Figure 1) and therefore the above written associativity law does not satisfy.
\begin{center}
\includegraphics[scale=0.6]{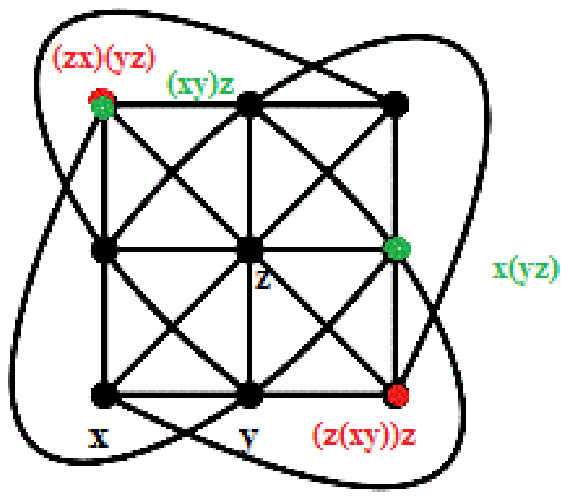}

{\bf Figure 1. {\it Points corresponding to the associativity $(xy)z=x(yz)$ are
marked in green; those  corresponding to the Moufang identity $(zx)(yz)=(z(xy))z$ in red.}}
\end{center}

Now, assume the opposite case. Here either (i) $x\neq y\neq z\neq x$ form a block or (ii) precisely two
among $x,y,z$ coincide or (iii) $x=y=z$. In case (i) one of elements $x, y, z$
is the product of the remaining two, say $z=xy$. In case (ii) there are two distinct
points among $x, y, z$, say $x\neq y=z$. Therefore, in cases (i), (ii) elements $\{x, y, xy, e\}$ form a
group isomorphic to the Klein four-group.

So, in cases (i), (ii),
$$[(x,\xi)\cdot (y,\eta)]\cdot (xy,\zeta) = (x,\xi)\cdot [(y,\eta)\cdot (xy,\zeta)]$$
for all $x, y, xy\in \{x, y, xy, e\}$ and $\xi, \eta, \zeta \in Z_2$ if and only if $f(x,x)=-1$ for all $x\in S\setminus\{e\}$.
That is, in cases (i), (ii)
$$<(a,\alpha); a\in \{x, y\}, \alpha \in Z_2>$$
forms a group isomorphic to the quaternion group precisely if $f(x,x)=-1$ for all $x\in S\setminus\{e\}$
(cf. Proposition 2 (b) in \cite{SS2}, p. 135). This means that in cases (i), (ii) discussed above loop $L$ has exponent $4$.

Next, set $f(u,u)=1$, $u\in S$, holds in $L$, i.e., $L$ is of exponent $2$.
Then however all elements of $L$ associate in any order except elements $(a, \xi), (a,\eta), (c, \zeta)$ and elements $(a, \xi), (c,\eta), (c,\zeta)$
$a\neq c$, $a, c \in \{x, y, xy\}$, ($\xi, \eta, \zeta \in Z_2$) of $L$ in the given orders they do not generate a group (cf. Proposition 2 (a) in \cite{SS2}, p. 135). If they would generate a group it would be a noncommutative group of order $8$ such that all elements distinct from the identity element have order $2$ and such group does not exist. The only group of order $8$ with the given order statistic is the elementary abelian group $E_8$ (cf. \cite{ATLAS})

This completes the
analysis in cases (i), (ii).

Finally, in case (iii), $[(x,\xi)\cdot (x,\eta)]\cdot (x,\zeta)=(x,\xi)\cdot [(x,\eta)\cdot (x,\zeta)]$, for any $x\in S$
and $\xi, \eta, \zeta \in Z_2$; this holds since $<(x,\xi); x\neq e, \xi \in Z_2>$ is again isomorphic to the
Klein four-group for both exponents of $L$.

Summarizing the above facts, in an OHL of exponent $4$ only elements of the set
$\{(x, \pm 1), (y, \pm 1), (e, \pm 1), (xy, \pm 1)\}$ associate, which form a subgroup of $\,L\,$ isomorphic to
the quaternion group.
\ep

As a byproduct we showed

\bc
OHLs of exponent $2$ do not satisfy the Moufang's theorem.
\ec

\br
OHLs of exponent $4$ are flexible, alternative and have the inverse property. OHLs of exponent $2$ are
flexible and have the automorphic inverse and the cross inverse properties. (cf. Theorem $5$ in \cite{SS2}, p. $136$).
\er

\noindent
{\bf Acknowledgement}
\vskip 10pt
\noindent
The author has been supported by FAPESP Grant - process No 11/51845-5,
and expresses the deep gratitude to IMS, University of S\~{a}o Paulo,
Brazil, for the warm hospitality.

\vskip 5pt
\noindent Izabella Stuhl\\
University of S\~{a}o Paulo, 05508-090 S\~{a}o Paulo, SP, Brazil\\
University of Debrecen, H-4010 Debrecen, Hungary\\
{\it E-mail}: {\it {}izabella@ime.usp.br}\\

\begin{thebibliography}{30}

\bibitem{Br}
R. H. Bruck, A survey of binary systems. \emph{Ergebnisse der Mathematik und ihrer Grenzgebiete. Neue Folge, Heft 20. Reihe: Gruppentheorie
Springer-Verlag 1958.}

\bibitem{Bu}
F. Buekenhout, Une caracterisation des espaces affins bas\'ee sur la notion de droite.
\emph{Math. Z., 111 (1969), 367-371.}

\bibitem{CLRS}
C. J. Colbourn, M. L. Merlini Giuliani, A. Rosa, I. Stuhl,
Steiner loops and Moufang's theorem. \emph{(submitted)}.

\bibitem{ATLAS}
J. H. Conway, R. T. Curtis, S. P. Norton, R. A. Parker, R. A. Wilson,
ATLAS of finite groups,
\emph{Oxford University Press, (2003)}

\bibitem{DP}
J. W. Di Paola, When is a totally symmetric loop a group?
\emph{The Amer. Math. Monthly, 76 (1969), 249-252.}

\bibitem{H}
M. Hall, Automorphisms of Steiner triple systems. \emph{IBM J. Res. Develop., (1960), 460-472.}

\bibitem{L}
M. L. Merlini Giuliani, G. Souza dos Anjos, Steiner Loops Satisfying Moufang's Theorem. \emph{(submitted)}

\bibitem{Mo}
R. Moufang, Zur Struktur von Alternativk\"orpern. \emph{Math. Ann. 110 (1935), 416-430.}

\bibitem{NS2}
P. T. Nagy, K. Strambach, Schreier loops.
\emph{Czechoslovak Math. J., 58 (133), (2008) 759-786.}

\bibitem{NSt2} P. T. Nagy, I. Stuhl,
Quasigroups arisen by right nuclear extension.
\emph{Comment. Math. Univ. Carolin., 53, (2012) 391-395.}

\bibitem{NSt} P. T. Nagy, I. Stuhl, Right nuclei of quasigroup extensions.
\emph{Comm. in Alg., 40, (2012) 1893-1900.}

\bibitem{RR}
R. Roth, D. K. Ray-Chaudhuri, Hall triple systems and commutative Moufang exponent $3$ loops: the case of nilpotence class $2$.
\emph{J. of Combinatorial Theory, Ser. A, 36, (1984), 129-162.}

\bibitem{SS2}
K. Strambach, I. Stuhl, Oriented Steiner loops.
\emph{Beitr\"age zur Algebra und Geometrie, 54, (2013) 131-145.}

\bibitem{SS}
K. Strambach, I. Stuhl, Translation groups of Steiner loops.
\emph{Discrete Mathematics, 309, (2009) 4225 -- 4227.}

\bibitem{S2}
I. Stuhl, Oriented Steiner quasigroups.
\emph{J. of Algebra and its Applications, 13, (2014)}.

\bibitem{S}
I. Stuhl, Regular permutations of quasigroups.
\emph{J. of Adv. Math. Studies, 3, (2010) 111-116.}

\bibitem{W}
S. Walter, Geordnete Steinersche Tripelsysteme. \emph{
Dissertation, Universit\"at Erlangen-N\"urnberg (1983).}

\end{thebibliography}
\end{document}